\documentclass[secnum,leqno,12pt]{article}
\usepackage[dvips]{graphicx}
\usepackage{latexsym}
\usepackage{amssymb}
\usepackage{amsmath}
\usepackage{amscd}
\usepackage{mathrsfs}
\usepackage{bm}
\usepackage{eepic}
\usepackage{epic}
\newtheorem{Th}{Theorem}[section]

\newtheorem{Lem}[Th]{Lemma}

\newtheorem{Rem}[Th]{Remark}

\newtheorem{Pro}[Th]{Proposition}
\newcommand{\demo}{\par\noindent{\it Proof. \/}\ }
\newcommand{\enD}{\hfill $\Box$ \vspace{3truemm}\par}
\newcommand{\bx}{\mbox{\boldmath $x$}}
\newcommand{\bX}{\mbox{\boldmath $X$}}
\newcommand{\be}{\mbox{\boldmath $e$}}

\newcommand{\bv}{\mbox{\boldmath $v$}}
\newcommand{\sbv}{\mbox{\scriptsize \boldmath$v$}}

\newcommand{\sbxi}{\mbox{\scriptsize \boldmath$\xi$}}
\newcommand{\by}{\mbox{\boldmath $y$}}
\newcommand{\bo}{\mbox{\boldmath $0$}}

\newcommand{\bw}{\mbox{\boldmath $w$}}

\newcommand{\bxi}{\mbox{\boldmath $\xi$}}

\newcommand{\bn}{\mbox{\boldmath $n$}}

\newcommand{\R}{{\mathbb R}}
\newcommand{\lon}{\longrightarrow}
\newcommand{\ou}{{\overline{u}}}

\title{ Geometry of world sheets in Lorentz-Minkowski space}
\author{ Shyuichi IZUMIYA \footnote{Department of Mathematics, Hokkaido University, Sapporo
060-0810, Japan.\newline e-mail: \texttt{izumiya@math.sci.hokudai.ac.jp}
}}
\begin{document}
\maketitle

\begin{abstract}
A world sheet in Lorentz-Minkowski space is a timelike submanifold consisting of a one-parameter family of spacelike submanifolds
in Lorentz-Minkowski space.
In this paper we investigate differential geometry of world sheets in Lorentz-Minkowski
space as an application of the theory of big wave fronts.
\end{abstract}

\section{Introduction}
\par
In this paper we consider differential geometry of world sheets in
Lorentz-Minkowski space. 
A world sheet is a timelike submanifold consisting of a one-parameter family of spacelike submanifolds
in a Lorentz manifold.
Since we do not have the notion of constant time in the relativity theory,
we  consider one-parameter families of spacelike submanifolds depending on the time-parameter (i.e., world sheets).
In this case, the spacelike submanifold with the constant parameter is not necessarily the constant time
in the ambient space. If we observe a surface in our space, then it is moving around the sun. Moreover, the solar system itself is moving depending on the Galaxy movement. Therefore, even if it looks a fixed surface (for example, a surface of a solid body) in Euclidean $3$-space, it is a three dimensional world sheets in
Lorentz-Minkowski $4$-space. Moreover, there appeared higher dimensional Lorentz manifolds in the theoretical physics (i.e., the super string theory, the brane world scenario etc.).
So we consider world sheets with general codimension in general dimensional Lorentz-Minkowski space.
In \cite{Izu-Kase} lightlike flat geometry on a spacelike submanifold with general codimension has been investigated.
Their method is quite useful for the study of the geometry of world sheets.
\par
On the other hand,
Lorentz-Minkowski space gives a geometric framework of the special relativity theory. Although there are no gravity
in Lorentz-Minkowski space, it provides a simple model of general Lorentz manifolds.
In this paper we investigate the lightlike geometry of world sheets in Lorentz-Minkowski space
with general codimension from the view point of the contact with lightlike hyperplanes.
The natural connection between geometry and singularities relies on the basic fact that the contact of a
submanifold with the models of the ambient space can be described by means of the analysis
if the singularities of appropriate families of contact functions, or equivalently, of their associated 
Lagrangian/Legendrian maps.
For the lightlike geometry the models are lightlike hyperplanes or lightcones.
The lightlike flat geometry is the lightlike geometry which adopts lightlike hyperplanes as model hypersurfaces.
Since we consider world sheets (i.e., one parameter families of spacelike submanifolds), 
the models are families of lightlike hyperplanes and 
the theory of one parameter bifurcations of Legendrian singularities is essentially useful.
Such a theory was initiated by Zakalyukin \cite{Zak76,Zak1} as the theory of big wave fronts.
There have been some developments on this theory during past two decades\cite{Izumiya93,Izu95,Izumiya-Takahashi,Izumiya-Takahashi2,Izumiya-Takahashi3, Zakalyukin95,Zakalyukin05}. Several applications of the theory were discovered in those articles.
For applying this theory, some equivalence relations among big wave fronts were used.
Here, we consider another equivalence relation among big wave fronts which is different from
the equivalence relations considered in those articles. This equivalence relation is corresponding to
the equivalence relation introduced in \cite{Go-Sch,Izudc} for applying the singularity theory to bifurcation problems.
\par
In \S 2 basic notations and properties of Lorentz-Minkowski space are explained.
Differential geometry of world sheets in Lorentz-Minkowski space is constructed in
\S 3.
We introduce the notion of (world and momentary) lightcone Gauss maps and induce the corresponding curvatures of world sheets respectively.
In \S 4 we define the lightcone height functions family and the extended lightcone height functions family of
a world sheet. We calculate the singular points of these families of functions and induce the notion of lightcone pedal maps and
unfolded lightcone pedal maps respectively. We investigate the geometric meanings of the singular points of the lightcone pedal maps
from the view point of the contact with families of lightlike hyperplanes in \S 5.
We can show that the image of the unfolded lightcone pedal map is a big wave front of a certain big Legendrian submanifold.
Therefore, we apply the theory of big wave fronts to our situation and interpret the geometric meanings of 
the singularities of the unfolded lightcone pedal map in \S 6.

\section{Basic concepts}
\par 
We introduce in this section some basic notions on Lorentz-Minkowski
$(n+1)$-space. For basic
concepts and properties, see \cite{Oneil}.
Let $\R^{n+1}=\{(x_0,x_1,\dots  ,x_n)\ |\ x_i \in \R\ (i=0,1,\dots , n)\}$
be an $(n+1)$-dimensional cartesian space. For any
$\bx =(x_0,x_1,\dots ,x_n),\ \by =(y_0,y_1,\dots ,y_n)\in \R^{n+1},$
the {\it pseudo scalar product } of $\bx $ and $\by $ is defined by
$
\langle \bx ,\by\rangle =-x_0y_0+\sum _{i=1}^n x_iy_i.
$
We call $(\R^{n+1} ,\langle ,\rangle )$  {\it Lorentz-Minkowski
$(n+1)$-space}. We write $\R^{n+1}_1$ instead of $(\R^{n+1} ,\langle
,\rangle )$. We say that a non-zero vector $\bx \in \R^{n+1}_1$ is
{\it spacelike, lightlike or timelike} if $\langle \bx,\bx \rangle
>0,$ $\langle \bx,\bx \rangle =0$ or  $\langle \bx,\bx \rangle <0$
respectively. The norm of the vector $\bx \in \R^{n+1}_1$ is defined
to be $\|\bx\|=\sqrt{|\langle \bx,\bx\rangle |}.$ We have the canonical
projection $\pi :\R^{n+1}_1\lon \R^n$ defined by $\pi (x_0,x_1,\dots
,x_n)=(x_1,\dots ,x_n).$ Here we identify $\{\bo\}\times \R^n$ with
$\R^n$ and it is considered as Euclidean $n$-space whose scalar
product is induced from the pseudo scalar product $\langle ,\rangle
.$ For a non-zero vector $\bv \in \R^{n+1}_1$ and a real number $c,$ we
define a {\it hyperplane with pseudo normal\/} $\bv$ by
$$
HP(\bv ,c)=\{\bx\in \R^{n+1}_1\ |\ \langle \bx ,\bv \rangle = c\ \}.
$$
We call $HP(\bv ,c)$ a {\it spacelike hyperplane\/}, a {\it timelike hyperplane\/}
or a {\it lightlike hyperplane\/}  if $\bv $ is timelike, spacelike or lightlike respectively.
\par
We now define {\it Hyperbolic $n$-space} by
$$
H^n_+(-1)=\{\bx\in \R^{n+1}_1 | \langle \bx ,\bx\rangle =-1, x_0>0
\}
$$
and {\it de Sitter $n$-space} by
$$
S^n_1=\{\bx\in \R^{n+1}_1 | \langle \bx ,\bx\rangle =1\ \}.
$$

We define
$$
LC^*=\{\bx=(x_0,x_1,\dots ,x_n)\in \R^{n+1}_1  \ | x_0\not= 0,\
\langle\bx ,\bx \rangle =0\}
$$
and we call it {\it the {\rm (}open{\rm )} lightcone} at the origin.
In the lightcone, we have the canonical unit spacelike sphere
defined by
\[
S^{n-1}_+=\{\bx =(x_0,x_1,\dots ,x_n)\ |\
\langle \bx ,\bx \rangle =0,\ x_0=1\}.
\]
We call $S^{n-1}_+$ the {\it lightcone unit $(n-1)$-sphere.}
If $\bx =(x_0,x_1,\dots ,x_n)$ is a lightlike vector, then
$x_0\not= 0.$ Therefore we have
\[
\widetilde{\bx}=\left(1,\frac{x_1}{x_0},\dots
,\frac{x_n}{x_0}\right)\in S^{n-1}_+.
\]
It follows that we have a projection
$\pi ^L_S:LC^*\lon S^{n-1}_+$ defined by $\pi ^L_S(\bm{x})=\widetilde{\bm{x}}.$
\par
For any $\bx_1,\bx_2,\dots ,\bx_n \in \R^{n+1}_1,$
we define a vector $\bx_1\wedge\bx_2\wedge\dots \wedge\bx_n$
by
\[
\bx_1\wedge\bx_2\wedge\dots \wedge\bx_n=
\left|
\begin{array}{cccc}
-\be_{0}&\be_{1}&\cdots &\be_{n}\\
x^1_{0}&x^1_{1}&\cdots &x^1_{n}\\
x^2_{0}&x^2_{1}&\cdots &x^2_{n}\\
\vdots &\vdots &\cdots &\vdots \\
x^n_{0}&x^n_{1}&\cdots &x^n_{n}
\end{array}
\right| ,
\]
where $\be_{0},\be_{1},\dots ,\be_{n}$ is the canonical basis of $\R^{n+1}_1$
and $\bx_i=(x_0^i,x_1^i,\dots ,x_n^i).$
We can easily check that
$
\langle \bx,\bx_1\wedge\bx_2\wedge\dots \wedge\bx_n\rangle ={\rm det}(\bx,\bx_1,\dots ,\bx_n),
$
so that
$\bx_1\wedge\bx_2\wedge\dots \wedge\bx_n$ is pseudo orthogonal to  any $\bx _i$ $(i=1,\dots ,n).$

\section{World sheets in Lorentz-Minkowski space}
\par
In this section we introduce the basic geometrical framework for the 
study of world sheets in Lorentz-Minkowski $(n+1)$-space.
Let $\R^{n+1}_1$ be a time-oriented 
space (cf., \cite{Oneil}). We choose $\be _0=(1,0,\dots ,0)$ as the future timelike
vector field. 
The world sheet is defined to be a timelike submanifold foliated by
a codimension one spacelike submanifolds.
Here, we only consider the local situation, so that we considered a one-parameter family of spacelike submanifolds.
Let $\bX :U\times I\lon \R^{n+1}_1$ be a timelike embedding of codimension $k-1,$
where $U\subset \R^s$ ($s+k=n+1$) is an open subset and $I$ an open interval.
 We write 
$W=\bX(U\times I)
$ and identify $W$ and $U\times I$ through the embedding $\bX.$
The embedding $\bX$ is said to be {\it timelike} if the tangent space $T_p  W$
of $W$ is a timelike subspace (i.e., Lorentz subspace of $T_p\R_1^{n+1}$) at any point $p\in W$.
We write $\mathcal{S}_t=\bX(U\times\{t\})$ for each $t\in I.$
We have a foliation $\mathcal{S}=\{\mathcal{S}_t\ |t\in I\}$
on $W$.
We say that $\mathcal{S}_t$ is {\it spacelike} if the tangent space $T_p\mathcal{S}_t$ 
consists only spacelike vectors (i.e., spacelike subspace) for
any point $p\in \mathcal{S}_t.$
We say that $(W, \mathcal{S})$ (or, $\bX$ ) is a {\it world sheet}
if $W$ is time-orientable and each $\mathcal{S}_t$ is spacelike.
We call $\mathcal{S}_t$ a {\it momentary space} of $(W,\mathcal{S})$.
For any $p=\bX(\overline{u},t)\in W\subset \R_1^{n+1},$
we have
\[
T_pW=\langle \bX_t(\ou,t),\bX_{u_1}(\ou,t),\dots ,\bX_{u_s}(\ou,t)\rangle _\R,
\]
where we write $(\overline{u},t)=(u_1,\dots ,u_s,t)\in U\times I$, $\bX_t=\partial \bX/\partial t$
and $\bX_{u_j}=\partial \bX/\partial u_j.$
We also have
\[
T_p\mathcal{S}_t=\langle \bX_{u_1}(\ou,t),\dots ,\bX_{u_s}(\ou,t)\rangle _\R.
\]
Since $W$ is time-orientable, there exists a timelike vector field $\bv(\ou,t)$ on $W$ \cite[Lemma 32]{Oneil}.
Moreover, we can choose that $\bv$ is {\it future directed} which means that
$\langle \bm{v}(\ou,t),\bm{e}_0\rangle <0.$
\par
Let $N_p(W)$ be the pseudo-normal space of $W$ at $p=\bX(\ou,t)$ in $\R^{n+1}_1.$
Since $T_pW$ is a timelike subspace of $T_p\R^{n+1}_1,$
$N_p(W)$ is a $(k-1)$-dimensional spacelike subspace of $T_p\R^{n+1}_1$
(cf.,\cite{Oneil}). On the pseudo-normal space $N_p(W),$ we have a $(k-2)$-unit sphere
\[
N_1(W)_p=\{\bxi\in N_p(W)\ |\ \langle \bxi,\bxi\rangle =1\ \}.
\]
Therefore, we have a unit spherical normal bundle over $W$:
\[
N_1(W)=\bigcup _{p\in W} N_1(W)_p.
\]
\par
On the other hand,
we write $N_p(\mathcal{S}_t)$ as the pseudo-normal space of $\mathcal{S}_t$ at $p=\bX(u,t)$
in $\R^{n+1}_1.$
Then $N_p(\mathcal{S}_t)$ is a $k$-dimensional Lorentz subspace of $T_p\R^{n+1}_1$.
On the pseudo-normal space $N_p(\mathcal{S}_t),$ we have two kinds of pseudo spheres:
\begin{eqnarray*}
N_p(\mathcal{S}_t;-1)& = & \{\bv\in N_p(\mathcal{S}_t)\ |\ \langle \bv,\bv\rangle =-1\ \} \\
N_p(\mathcal{S}_t;1)&= & \{\bv\in N_p(\mathcal{S}_t)\ |\ \langle \bv,\bv\rangle =1\ \}.
\end{eqnarray*}
We remark that $N_p(\mathcal{S}_t;-1)$ is the $(k-1)$-dimensional hyperbolic space and $N_p(\mathcal{S}_t;1)$ is the $(k-1)$-dimensional de Sitter space.
Therefore,
we have two unit spherical normal bundles $N(\mathcal{S}_t;-1)$ and $N(\mathcal{S}_t;1)$ over $\mathcal{S}_t$. 
Since $\mathcal{S}_t=\bX(U\times \{t\})$ is a codimension one spacelike submanifold in $W,$ there exists a
unique timelike future directed unit normal vector field $\bn^T(\ou,t)$ of $\mathcal{S}_t$ such that 
$\bn^T(\ou,t)$ is tangent to $W$ at any point $p=\bX(\ou,t).$
It means that $\bn^T(\ou,t)\in N_p(\mathcal{S}_t)\cap T_pW$ with $\langle \bn^T(\ou,t),\bn^T(\ou,t)\rangle =-1$ and
$\langle\bn^T(\ou,t),\be _0\rangle <0.$
We define a $(k-2)$-dimensional spacelike unit sphere in $N_p(\mathcal{S}_t)$ by
\[
N_1(\mathcal{S}_t)_p[\bn ^T]=\{\bxi \in N_p(\mathcal{S}_t;1)\ |\ \langle \bxi, \bn ^T(\ou,t)\rangle =0,p=\bX(\ou,t)\ \}.
\]
Then we have a {\it spacelike unit $(k-2)$-spherical bundle $N_1(\mathcal{S}_t)[\bn ^T]$ over $\mathcal{S}_t$ with respect to $\bn ^T$}.
Since we have
$T_{(p,\xi)}N_1(\mathcal{S}_t)[\bn^T]=T_p\mathcal{S}_t\times T_\xi N_1(\mathcal{S}_t)_p[\bn ^T],$
we have the canonical Riemannian metric on $N_1(\mathcal{S}_t)[\bn^T]$
which we write $(G_{ij}((\ou,t),\bxi))_{1\leqslant i,j\leqslant n-1}.$
Since $\bn^T$ is uniquely determined, we write $N_1[\mathcal{S}_t]=N_1(\mathcal{S}_t)[\bn^T].$
Moreover, we remark that $N_1(W)|\mathcal{S}_t=N_1[\mathcal{S}_t]$ for any $t\in I.$

\par
We now define a map $\mathbb{LG}:N_1(W)\lon LC^*$ by $\mathbb{LG}(\bX(\ou,t),\bxi)=\bn^T(\ou,t)+\bxi$.
We call $\mathbb{LG}$ a {\it world lightcone Gauss map} of
$N_1(W)$, where $W=\bX(U\times I)$.
A {\it momentary lightcone Gauss map} of $N_1[\mathcal{S}_t]$ is defined 
to be the restriction of the world lightcone Gauss map of $N_1(W)$:
\[
\mathbb{LG}(\mathcal{S}_t)=\mathbb{LG}|N_1[\mathcal{S}_t]:N_1[\mathcal{S}_t]\lon LC ^*.
\]
This map leads us to the notions of curvatures.
Let $T_{(p,\bm{\xi})}N_1[\mathcal{S}_t]$ be the tangent space of $N_1[\mathcal{S}_t]$ at $(p,\bxi).$
With the canonical identification 
$$
(\mathbb{LG}(\mathcal{S}_t)^*T\R^{n+1}_1)_{(p,\bm{\xi})}
=T_{(\bm{n}^T(p)+\bm{\xi})}\R^{n+1}_1\equiv T_p\R^{n+1}_1,
$$
we have
\[
T_{(p,\sbxi)}N_1[\mathcal{S}_t]=T_p\mathcal{S}_t\oplus T_\xi S^{k-2}\subset T_p\mathcal{S}_t\oplus N_p(\mathcal{S}_t)=T_p\R^{n+1}_1,
\]  
where $T_\xi S^{k-2}\subset T_\xi N_p(\mathcal{S}_t)\equiv N_p(\mathcal{S}_t)$ and $p=\bX(\overline{u},t).$
Let 
\[
\Pi ^t :\mathbb{LG}(\mathcal{S}_t)^*T\R^{n+1}_1=TN_1[\mathcal{S}_t]\oplus \R^{s+2}
\lon TN_1[\mathcal{S}_t]
\]
be the canonical projection.
Then
we have
a linear transformation
\[
S_\ell (\mathcal{S}_t)_{(p,\sbxi)}=-\Pi^t_{\mathbb{LG}(\mathcal{S}_t)(p,\xi)}\circ d_{(p,\xi)}\mathbb{LG}(\mathcal{S}_t)
: T_{(p,\xi)}N_1[\mathcal{S}_t]\lon T_{(p,\xi)}N_1[\mathcal{S}_t],
\]
which is called the {\it momentary lightcone shape operator} of $N_1[\mathcal{S}_t]$ at $(p,\bxi).$ 
\par
On the other hand,  for $t_0\in I,$ we choose a spacelike unit vector field $\bm{n}^S$ along $W=\bm{X}(U\times I)$ at
least locally such that 
$\bn ^S(u,t_0)\in N_1(\mathcal{S}_{t_0}).$ 
Then we have 
$\langle \bn ^S,\bn ^S\rangle =1$ and $\langle \bX_t,\bn^S\rangle = \langle \bX_{u_i},\bn^S\rangle =\langle \bn ^T,\bn ^S\rangle =0$ at $(\ou,t_0)\in U\times I.$
Clearly, the vector
$\bn^T (\ou,t_0)+ \bn^S(\ou,t_0)$ is lightlike. 
We define a mapping
\[
\mathbb{LG}(\mathcal{S}_{t_0};\bn^S):U\lon LC^*
\]
by $\mathbb{LG}(\mathcal{S}_{t_0};\bn^S)(u)=\bn^T(\ou,t_0)+\bn^S(\ou,t_0),$
which is called a {\it momentary lightcone Gauss map of
$\mathcal{S}_{t_0}=\bX(U\times \{t_0\})$ with respect to
$\bn^S.$}
Under the identification of $\mathcal{S}_{t_0}$ and $U\times\{t_0\}$ through $\bX,$ we have the
linear mapping provided by the derivative of the momentary lightcone Gauss map $\mathbb{LG}(\mathcal{S}_{t_0};\bn^S)$ at each point 
$p=\bX(\ou,t_0)$,
\[
d_p\mathbb{LG}(\mathcal{S}_{t_0};\bn^S):T_p\mathcal{S}_{t_0}\lon T_p\R^{n+1}_1= T_p\mathcal{S}_{t_0}\oplus N_p(\mathcal{S}_{t_0}).
\]
Consider the orthogonal projection $\pi ^t:T_p\mathcal{S}_{t_0}\oplus
N_p(\mathcal{S}_{t_0})\rightarrow T_p\mathcal{S}_{t_0}.$ We define
\[
S_p(\mathcal{S}_{t_0};\bn^S)=-\pi^t\circ d_p\mathbb{LG}(\mathcal{S}_{t_0};\bn^S):T_p\mathcal{S}_{t_0}\lon T_p\mathcal{S}_{t_0}.
\]
We call the
linear transformation $S_{p}(\mathcal{S}_{t_0};\bn^S)$ an {\it $\bn^S$-momentary shape
operator} of $\mathcal{S}_{t_0}=\bX (U\times \{t_0\})$ at $p=\bX (\overline{u},t_0).$ 
Let $\{\kappa
_{i}(\mathcal{S}_{t_0};\bn^S)(p)\}_ {i=1}^s$ be the eigenvalues of $S_{p}(\mathcal{S}_{t_0};\bn^S)$, which are called {\it momentary lightcone
principal curvatures of $\mathcal{S}_{t_0}$ with respect to $\bn^S$\/} at $p=\bX(\overline{u},t_0)$.
Then a {\it momentary lightcone Lipschitz-Killing curvature of $\mathcal{S}_{t_0}$ with respect to
$\bn^S$\/} at $p=\bX (\overline{u},t_0)$ is defined as follows:
\[
K_\ell(\mathcal{S}_{t_0};\bn^S)(p)={\rm det} S_{p}(\mathcal{S}_{t_0};\bn^S).
\]
We say that a point $p=\bX (\ou,t_0)$ is an {\it $\bn^S$-momentary lightcone umbilical
point} of $\mathcal{S}_{t_0}$ if 
\[
S_{p}(\mathcal{S}_{t_0};\bn^S)=\kappa (\mathcal{S}_{t_0};\bn^S)(p) 1_{T_{p}\mathcal{S}_{t_0}}.
\]
We say that $W=\bX (U\times I)$ is {\it totally
$\bn^S$-lightcone umbilical} if each point $p=\bX(\ou,t)\in W$ is
an $\bn^S$-momentary lightcone umbilical point of $\mathcal{S}_t .$
Moreover, $W=\bX(U\times I)$ is said to be {\it totally lightcone umbilical} if
it is totally $\bn^S$-lightcone umbilical for any $\bn^S.$
We deduce now the lightcone Weingarten formula. Since $\mathcal{S}_{t_0}=\bX(U\times \{t_0\})$
is spacelike submanifold, we have a Riemannian metric
(the {\it first fundamental form \/}) on $\mathcal{S}_{t_0}$
defined by $ds^2 =\sum _{i=1}^{s} g_{ij}du_idu_j$,  where
$g_{ij}(\ou,t_0) =\langle \bX _{u_i}(\ou,t_0 ),\bX _{u_j}(\ou,t_0)\rangle$ for any
$\ou\in U.$ We also have a {\it lightcone second fundamental invariant of $\mathcal{S}_{t_0}$
with respect to the normal vector field $\bn ^S $\/} defined
by $h _{ij}(\mathcal{S}_{t_0};\bn^S )(\ou,t_0)=\langle -(\bn^T +\bn^S)
_{u_i}(\ou,t_0),\bX_{u_j}(\ou,t_0)\rangle$ for any $\ou\in U.$
By the similar arguments to those in the proof of \cite[Proposition 3.2]{IzuSM}, we have 
the following proposition.
\begin{Pro}
We choose a pseudo-orthonormal frame $\{\bn^T,\bn^S_1,\dots ,\bn^S_{k-1}\}$ of $N(\mathcal{S}_{t_0})$ with $\bn^S_{k-1}=\bn^S.$ Then we have the following lightcone Weingarten formula {\rm :}
\vskip1.5pt
\par\noindent
{\rm (a)} $\mathbb{LG}(\mathcal{S}_{t_0};\bn^S)_{u_i}=\langle \bn ^T_{u_i},\bn ^S\rangle(\bn^T+\bn^S)+\sum _{\ell =1}^{k-2}\langle (\bn^T+\bn^S)_{u_i},\bn^S_\ell \rangle\bn^S_\ell  \\ \hskip2.8cm-\sum_{j=1}^{s}
h_i^j(\mathcal{S}_{t_0};\bn^S )\bX _{u_j}$,
\par\noindent
{\rm (b)} $
\pi ^t\circ \mathbb{LG}(\mathcal{S}_{t_0};\bn^S)_{u_i}=-\sum_{j=1}^{s}
h_i^j(\mathcal{S}_{t_0};\bn^S )\bX _{u_j}.
$
\smallskip
\par\noindent
Here, $\displaystyle{\left(h_i^j(\mathcal{S}_{t_0};\bn^S )\right)=\left(h_{ik}(\mathcal{S}_{t_0};\bn^S)\right)\left(g^{kj}\right)}$
and $\displaystyle{\left( g^{kj}\right)=\left(g_{kj}\right)^{-1}}.$
\end{Pro}
\par
Since $\mathbb{LG}(\mathcal{S}_{t_0};\bn^S)_{u_i}=d\mathbb{LG}(\mathcal{S}_{t_0};\bn^S)(\bX_{u_i})$,
we have 
\[
S_p(\mathcal{S}_{t_0};\bn^S )(\bX_{u_i}(\ou,t_0))=-\pi^t\circ \mathbb{LG}(\mathcal{S}_{t_0};\bn^S)_{u_i}(\ou,t_0),
\]
so that the representation matrix of $S_p(\mathcal{S}_{t_0};\bn^S )$ with respect to the basis 
$
\{\bX_{u_i}(\ou,t_0)\}_{i=1}^s$ of $T_p\mathcal{S}_{t_0}$ is $(h^i_j(\mathcal{S}_{t_0};\bn^S)(\ou,t_0)).$
Therefore, we have an explicit
expression of the momentary lightcone Lipschitz-Killing curvature of $\mathcal{S}_{t_0}$ with respect to
$\bn^S$ as follows:
$$
K_\ell (\mathcal{S}_{t_0};\bn^S )(\ou,t_0)=\frac{\displaystyle{{\rm det}\left(h_{ij}(\mathcal{S}_{t_0};\bn^S )(\ou,t_0)\right)}}
{\displaystyle{{\rm det}\left(g_{\alpha \beta}(\ou,t_0)\right)}}.
$$
Since $\langle -(\bn^T +\bn^S )(\ou,t_0),\bX _{u_j}(\ou,t_0)\rangle =0$,  we have
\[
h_{ij}(\mathcal{S}_{t_0};\bn^S)(\ou,t_0)=\langle \bn^T (\ou,t_0)+\bn^S (\ou,t_0),\bX
_{u_iu_j}(\ou,t_0)\rangle.
\]
 Therefore the lightcone second fundamental
invariants of $\mathcal{S}_{t_0}$ at a point $p_0=\bX (\ou_0,t_0)$ depend only on the values 
$\bn^T (\ou_0)+\bn^S (\ou_0)$ and $\bX _{u_iu_j}(\ou_0)$, respectively.
Therefore, we write
\[
h_{ij}(\mathcal{S}_{t_0};\bn^S)(\ou_0,t_0)=h_{ij}(\mathcal{S}_{t_0})(p_0,\bxi_0),
\]
where $p_0=\bX(\ou_0,t_0)$ and $\bxi_0=\bn^S(\ou_0,t_0)\in N_1(W)_{p_0}.$
Thus, the $\bn^S$-momentary shape operator and the momentary lightcone curvatures also depend only on
$\bn^T (\ou_0,t_0)+\bn^S (\ou_0,t_0)$, $\bX_{u_i}(\ou_0,t_0)$  and $\bX
_{u_iu_j}(\ou_0,t_0)$, which are independent of the derivations of the vector fields
$\bn^T$ and 
$\bn^S .$ 
It follows that we write $S_{p_0}(\mathcal{S}_{t_0};\bxi_0)=S_{p_0}(\mathcal{S}_{t_0};\bn^S),$ $\kappa _i(\mathcal{S}_{t_0},\bxi_0)(p_0)= \kappa _i(\mathcal{S}_{t_0};\bn^S)(p_0)$ $(i=1,\dots ,s)$
and $K_\ell(\mathcal{S}_{t_0},\bxi_0)(p_0)=K_\ell (\mathcal{S}_{t_0};\bn^S)(p_0)$ at $p_0=\bX (u_0,t_0)$
with respect to $\bxi_0=\bn^S(u_0,t_0).$ 
We also say that a point $p_0=\bX (u_0,t_0)$ is 
{\it $\bxi_0$-momentary lightcone umbilical\/} if $S_{p_0}(\mathcal{S}_{t_0};\bxi_0)=\kappa _i(\mathcal{S}_{t_0})(p_0,\bxi_0)1_{T_{p_0}\mathcal{S}_{t_0}}$. 
We say that a point $p_0=\bX (u_0,t_0)$ is a {\it $\bxi_0$-momentary lightcone parabolic point \/} of $\mathcal{S}_{t_0}$ if
$K_\ell (\mathcal{S}_{t_0};\bxi_0)(p_0)=0.$ 
\par
Let $\kappa _\ell(\mathcal{S}_{t})_i(p,\bxi)$ be the eigenvalues of the lightcone shape operator $S_\ell(\mathcal{S}_{t}) _{(p,\sbxi)}$, $(i=1,\dots ,n-1)$. 
We write $\kappa _\ell (\mathcal{S}_t)_i(p,\bxi)$, $(i=1,\dots ,s)$ for the eigenvalues whose 
eigenvectors belong to $T_p\mathcal{S}_{t}$
and $\kappa _\ell (\mathcal{S}_t)_i(p,\bxi)$, $(i=s+1,\dots n)$ for the eigenvalues whose eigenvectors belong to 
the tangent space of the fiber  
of $N_1[\mathcal{S}_t].$  
\begin{Pro}
For $p_0=\bX(u_0,t_0)$ and $\bxi_0\in N_1[\mathcal{S}_{t_0}]_{p_0},$ we have
$$\kappa _\ell (\mathcal{S}_{t_0})_i(p_0,\bxi_0)=\kappa _i(\mathcal{S}_{t_0},\bxi_0)(p_0),\ (i=1,\dots s),\ \kappa _\ell (\mathcal{S}_{t_0})_i(p_0,\bxi_0)=-1,\ (i=s+1,\dots n).$$
\end{Pro}
\demo
Since $\{\bn^T,\bn^S_1,\dots ,\bn^S_{k-1}\}$ is a pseudo-orthonormal frame of $N(\mathcal{S}_t)$
and $\bxi_0=\bn^S_{k-1}(\ou_0,t_0)\in S^{k-2}=N_1[\mathcal{S}_{t_0}]_p,$ 
we have $
\langle \bn^T(\ou_0,t_0),\bxi_0\rangle =\langle \bn^S_i(\ou_0,t_0),\bxi_0\rangle =0$
for $i=1,\dots ,k-2.$
Therefore, we have 
$$T_{\bm{\xi}_0}S^{k-2}=\langle \bn^S_1(\ou_0,t_0),\dots ,\bn^S_{k-2}(\ou_0,t_0)\rangle .$$
By this orthonormal basis of $T_{\sbxi_0}S^{k-2},$
the canonical Riemannian metric $G_{ij}(p_0,\bxi_0)$ is represented by
\[
(G_{ij}(p_0,\bxi))=\left(
\begin{array}{cc}
g_{ij}(p_0)  & 0 \\
0 & I_{k-2}
\end{array}
\right) ,
\]
where $g_{ij}(p_0)=\langle \bX_{u_i}(\ou_0,t_0), \bX_{u_j}(\ou_0,t_0)\rangle $.
\par
On the other hand, by Proposition 3.1, we have
\[
-\sum_{j=1}^s h^j_i(\mathcal{S}_{t_0},\bn^S)\bX_{u_j}=\mathbb{LG}(\mathcal{S}_{t_0},\bn^S)_{u_i}=
d_{p_0}\mathbb{LG}(\mathcal{S}_{t_0};\bn^S)\left(\frac{\partial}{\partial u_i}\right),
\]
so that we have
\[
S_\ell(\mathcal{S}_{t_0})_{(p_0,\bm{\xi}_0)}\left(\frac{\partial}{\partial u_i}\right)=\sum_{j=1}^s h^j_i(\mathcal{S}_{t_0},\bn^S)\bX_{u_j}.
\]
Therefore, the representation matrix of $S_\ell(\mathcal{S}_{t_0})_{(p_0,\bxi_0)}$ with respect to the basis
$$
\{\bX_{u_1}(\ou_0,t_0),\dots ,\bX_{u_s}(\ou_0,t_0),\bn^S_1(\ou_0,t_0),\dots ,\bn^S_{k-2}(\ou_0,t_0)\}
$$ of $T_{(p_0,\bm{\xi}_0)}N_1[\mathcal{S}_{t_0}]$
is of the form
\[
\left(
\begin{array}{cc}
h^j_i(\mathcal{S}_{t_0},\bn^S)(u_0,t_0)  & * \\
0 & -I_{k-2}
\end{array}
\right).
\]
Thus, the eigenvalues of this matrix are $\lambda _i=\kappa _i(\mathcal{S}_{t_0},\bxi_0)(p_0)$, $(i=1,\dots ,s)$ and
$\lambda _i=-1,$ $(i=s+1,\dots ,n-1)$.
This completes the proof.
\enD

We call $\kappa _\ell (\mathcal{S}_{t})_i(p,\bxi)=\kappa _i(\mathcal{S}_{t},\bxi)(p)$, $(i=1,\dots ,s)$
{\it momentary lightcone principal curvatures} of $\mathcal{S}_{t}$ with respect to $\bxi$ at $p=\bX(\overline{u},t)\in W.$
\par
On the other hand, we define a mapping $\widetilde{\mathbb{LG}}(\mathcal{S}_t):N_1(\mathcal{S}_t)\lon
S^{n-1}_+$ by 
\[
\widetilde{\mathbb{LG}}(\mathcal{S}_t)(p,\bm{\xi})=\pi^L_S( \mathbb{LG}(\mathcal{S}_t)(p,\bm{\xi})),
\]
 which is called a {\it normalized momentary lightcone Gauss map} of $N_1(\mathcal{S}_t).$
A {\it normalized momentary lightcone Gauss map of $\mathcal{S}_t$ with respect to $\bm{n}^S$} is a mapping
$\widetilde{\mathbb{LG}}(\mathcal{S}_t;\bm{n}^S):U\lon S^{n-1}_+$ 
defined to be
$\widetilde{\mathbb{LG}}(\mathcal{S}_t;\bm{n}^S)(\ou)=\pi^L_S(\mathbb{LG}(\mathcal{S}_t;\bm{n}^S)(\ou)).$
The normalized momentary lightcone Gauss map of $\mathcal{S}_t$ with respect to $\bm{n}^S$
also induces a linear mapping $d_p\widetilde{\mathbb{LG}}(\mathcal{S}_t;\bn^S):T_p\mathcal{S}_t\lon T_p\R^{n+1}_1$ under
the identification of $U\times \{t\}$ and $\mathcal{S}_t,$ where $p=\bX(\ou,t).$
We have the following proposition.
\begin{Pro}
With the above notations, we have the following normalized lightcone Weingarten formula with respect to $\bn^S$:
\[
\pi^t\circ \widetilde{\mathbb{LG}}(\mathcal{S}_t;\bn^S)_{u_i}(\ou)=-\sum_{j=1}^s \frac{1}{\ell _0(\ou,t)}h^j_i(\mathcal{S}_t;\bn^S)(\ou,t)\bX_{u_j}(\ou,t),
\]
where $\mathbb{LG}(\mathcal{S}_t;\bn^S)(\ou)=(\ell _0(\ou,t),\ell _1(\ou,t),\dots ,\ell _n\ou,t)).$
\end{Pro}
\demo
By definition, we have 
$\ell_0\widetilde{\mathbb{LG}}(\mathcal{S}_t;\bn^S)=\mathbb{LG}(\mathcal{S}_t;\bn^S).$ It follows that
\[
\ell _0\widetilde{\mathbb{LG}}(\mathcal{S}_t;\bn^S)_{u_i}=\mathbb{LG}(\mathcal{S}_t;\bn^S)_{u_i}-\ell _{0u_i}\widetilde{\mathbb{LG}}(\mathcal{S}_t;\bn^S).
\]
Since $\widetilde{\mathbb{LG}}(\mathcal{S}_t;\bn^S)(\ou)\in N_p(\mathcal{S}_t),$ we have
$
\pi^t\circ\widetilde{\mathbb{LG}}(\mathcal{S}_t;\bn^S)_{u_i}=\frac{1}{\ell _0}\pi^t\circ\mathbb{LG}(\mathcal{S}_t;\bn^S)_{u_i}.
$
By the lightcone Weingarten formula with respect to $\bn^S$ (Proposition 3.1), we have the desired formula.
\enD
\par
We call the linear transformation $\widetilde{S}_p(\mathcal{S}_t;\bm{n}^S)=-\pi^t\circ d_p\widetilde{\mathbb{LG}}(\mathcal{S}_t;\bn^S)$ a {\it normalized momentary
lightcone shape operator of $\mathcal{S}_t$ with respect to $\bn^S$} at $p$.
The eigenvalues $\{\widetilde{\kappa} _i(\mathcal{S}_t;\bn^S)(p)\}_{i=1}^s$ of $\widetilde{S}_p(\mathcal{S}_t;\bm{n}^S)$ are called {\it normalized momentary
lightcone principal curvatures}.
By the above proposition, we have
$\widetilde{\kappa} _i(\mathcal{S}_t;\bn^S)(p)=(1/\ell _0(\ou,t))\kappa _i(\mathcal{S}_t;\bn^S)(p).$
A {\it normalized momentary Lipschitz-Killing curvature} of $\mathcal{S}_t$ with respect to $\bn^S$ is defined to be
$\widetilde{K}_\ell (\ou,t)={\rm det}\, \widetilde{S}_p(\mathcal{S}_t;\bm{n}^S).$
Then we have the following relation between the normalized momentary lightcone Lipschitz-Killing curvature and the momentary lightcone
Lipschitz-Killing curvature:
\[
\widetilde{K}_\ell(\mathcal{S}_t;\bn^S)(p)=\left(\frac{1}{\ell_0(\ou,t)}\right)^sK_\ell (\mathcal{S}_t;\bn^S)(p),
\]
where $p=\bm{X}(\ou,t).$
By definition, $p_0=\bX(\ou_0,t_0)$ is the $\bn^S_0$-momentary umbilical point if and only if
$\widetilde{S}_{p_0}(\mathcal{S}_t;\bm{n}^S_0)=\widetilde{\kappa} _i(\mathcal{S}_{t_0};\bn^S)(p_0)1_{T_{p_0}\mathcal{S}_{t_0}}.$
We have the following proposition.
\begin{Pro} For any $t_0\in I,$
the following conditions {\rm (1)} and {\rm (2)} are equivalent{\rm :}
\par\noindent
{\rm (1)} There exists a spacelike unit vector field $\bn^S$ along $W=\bX(U\times I)$ 
such that $\bm{n}^S(\ou,t_0)\in N_1(\mathcal{S}_{t_0})$  and the normalized momentary lightcone Gauss map $\widetilde{\mathbb{LG}}(\mathcal{S}_{t_0};\bn^S)$ of $\mathcal{S}_{t_0}=\bX(U\times \{t_0\})$ with respect to
$\bn^S$ is constant. 
\par\noindent
{\rm (2)} There exists $\bv\in S^{n-1}_+$ and a real number $c$ such that $\mathcal{S}_{t_0}\subset HP(\bv,c).$
\par
Suppose that the above conditions hold. Then 
\par\noindent
{\rm (3)} $\mathcal{S}_{t_0}=\bX(U\times \{t_0\})$ is totally $\bn^S$-momentary flat.
\end{Pro}
\demo
Suppose that the condition (1) holds.
We consider a function $F:U\lon \R$ defined by $F(\ou)=\langle \bX(\ou,t_0),\bv\rangle.$
By definition, we have
\[
\frac{\partial F}{\partial u_i}(\ou)=\langle \bX_{u_i}(\ou,t_0),\bv\rangle
=\langle \bX_{u_i}(\ou,t_0),\widetilde{\mathbb{LG}}(\mathcal{S}_{t_0};\bn^S)(\ou)\rangle=0 ,
\]
for any $i=1,\dots ,s.$ Therefore, $F(\ou)=\langle \bX(\ou,t_0),\bv\rangle =c$ is constant.
It follows that $\mathcal{S}_{t_0}\subset HP(\bv, c)$ for $\bv\in S^{n-1}_+.$
\par
Suppose that $\mathcal{S}_{t_0}$ is a subset of a lightlike hyperplane $H(\bv,c)$ for $\bv\in S^{N-1}_+.$
Since $\mathcal{S}_{t_0}\subset HP(\bv ,c),$ we have $T_p\mathcal{S}_{t_0}\subset H(\bv ,0)$
for any $p\in \mathcal{S}_{t_0}.$
If $\langle \bn^T(\ou,t),\bv\rangle=0,$ then $\bn^T(\ou,t)\in HP(\bv ,0).$
We remark that $HP(\bv ,0)$ does not contain timelike vectors.
This is a contradiction. So we have $\langle \bn^T(\ou,t),\bv\rangle \not=0.$
We now define a vector field along $W=\bX(U\times I)$ by
\[
\bn^S(\ou,t)=\frac{-1}{\langle \bn^T(\ou,t),\bv\rangle}\bv-\bn^T(\ou,t).
\]
We can easily show that $\langle \bn^S(\ou,t),\bn^S(\ou,t)\rangle =1$ and $\langle \bn^S(\ou,t),\bn^T(\ou,t)\rangle =0.$
Since $T_p\mathcal{S}_{t_0}\subset H(\bv ,0)$, we have
$\langle \bm{X}_{u_i}(\ou,t_0),\bm{n}^S(\ou,t_0)\rangle =0.$
Hence $\bn^S$ is a spacelike unit vector field $\bn^S$ along $W=\bX(U\times I)$ 
such that $\bm{n}^S(\ou,t_0)\in N_1(\mathcal{S}_{t_0})$ and
$\widetilde{\mathbb{LG}}(\mathcal{S}_{t_0};\bn^S)(\ou)=\bv.$
By Proposition 3.3, if $\widetilde{\mathbb{LG}}(\mathcal{S}_{t_0};\bn^S)$ is constant,
then
$(h^j_i(\mathcal{S}_{t_0};\bn^S)(\ou,t_0))=O$. It follows that $\mathcal{S}_{t_0}$ is lightcone $\bn^S$-flat.
\enD

\section{Lightcone height functions}
In order to study the geometric meanings of the normalized lightcone Lipschitz-Killing curvature
$\widetilde{K}_\ell(\mathcal{S}_t;\bn^S)$ of $\mathcal{S}_t=\bX(U\times\{t\})$, we  introduce
a family of functions on $M=\bX(U).$
A family of {\it lightcone height functions} 
$
H:U\times (S^{n-1}_+\times I)\lon \R
$ on $W=\bX(U\times I)$ 
is defined to be
$
H((\ou,t),\bv)=\langle \bX(\ou,t),\bv\rangle.
$
The Hessian matrix of the lightcone height function $h_{(t_0,\bm{v}_0)}(\ou)=H((\ou,t_0),\bv_0)$
at $\ou_0$ is denoted by ${\rm Hess}(h_{(t_0,\bm{v}_0)})(\ou_0).$
The following proposition characterizes the lightlike parabolic points and lightlike flat points
in terms of the family of lightcone height functions.
\begin{Pro} Let  
$H:U\times (S^{n-1}_+\times I)\lon \R$ be the family of lightcone height functions on a world sheet $W=\bX(U\times I).$
Then 
\par\noindent
{\rm (1)} $(\partial H/\partial u_i)(\ou_0,t_0,\bv_0)=0\ (i=1,\dots ,s)$ if and only if
there exists a spacelike section $\bn^S$ of $N_1(\mathcal{S}_{t_0})$  such that $\bv_0=\widetilde{\mathbb{LG}}(\mathcal{S}_{t_0};\bn^S_0)(\ou_0).$ 
\par
Suppose that $p_0=\bX(\ou_0,t_0)$, $\bv_0=\widetilde{\mathbb{LG}}(\mathcal{S}_{t_0};\bn^S_0)(\ou_0).$ Then
\par\noindent
{\rm (2)} $p_0$ is an $\bn^S_0$-parabolic point of $\mathcal{S}_{t_0}$ if and only if
${\rm det}\,{\rm Hess}(h_{(t_0,\bm{v}_0)})\ou_0)=0,$ 
\par\noindent
{\rm (3)} $p_0$ is a flat $\bn^S_0$-umbilical point of $\mathcal{S}_{t_0}$ if and only if
${\rm rank}\,{\rm Hess}(h_{(t_0,\bm{v}_0)})\ou_0)=0$.
\end{Pro}
\demo
(1) Since $(\partial H/\partial u_i)((\ou_0,t_0)\bv_0)=\langle \bX_{u_i}(\ou_0,t_0),\bv_0),$
$(\partial H/\partial u_i)((\ou_0,t_0),\bv_0)=0\ (i=1,\dots ,s)$ if and only if $\bv_0\in N_{p_0}(\mathcal{S}_{t_0})$ and
$\bv_0\in S^{n-1}_+.$
By the same construction as in the proof of Proposition 3.4,
we have a spacelike unit normal vector field $\bn^S$ along $W=\bm{X}(U\times I)$ with
$\bm{n}^S(\ou,t_0)\in N_1(\mathcal{S}_{t_0})$ such that $\bv_0=\widetilde{\mathbb{LG}}(\mathcal{S}_{t_0};\bn^S)(\ou_0)
=\widetilde{\mathbb{LG}}(\mathcal{S}_{t_0};\bn^S_0)(\ou_0).$ 
The converse also holds.
For the proof of the assertions (2) and (3), as a consequence of Proposition 3.1, we have
\begin{eqnarray*}
{\rm Hess}(h_{(t_0,\sbv_0)})(\ou_0)&=&\left(\langle \bX_{u_iu_j}(\ou_0,t_0),\widetilde{\mathbb{LG}}(\mathcal{S}_{t_0};\bn^S)(\ou_0)\rangle\right) \\
&=&\left( \frac{1}{\ell _0}\langle \bX_{u_iu_j}(\ou_0,t_0),\bn^T(\ou_0,t_0)+\bn^S(\ou_0,t_0)\rangle\right) \\
&=& \left(\frac{1}{\ell _0}\langle \bX_{u_i}(\ou_0,t_0),(\bn^T+\bn^S)_{u_j}(\ou_0,t_0)\rangle\right) \\
&=&\left(\frac{1}{\ell _0}\langle \bX_{u_i}(\ou_0,t_0),-\sum _{k=1}^s h^k_j(\mathcal{S}_{t_0};\bn^S)(\ou_0)\bX_{u_k}(\ou_0,t_0)\rangle\right) \\
&=& \left(-\frac{1}{\ell _0}h_{ij}(\mathcal{S}_{t_0};\bn^S)(\ou_0)\right).
\end{eqnarray*}
By definition,  $K_\ell (\mathcal{S}_{t_0};\bn^S)(\ou_0)= 0$ if and only if ${\rm det}\, (h_{ij}(\mathcal{S}_{t_0};\bn^S)(\ou_0))=0.$
The assertion (2) holds.
Here, $p_0$ is a flat $\bn^S_0$-umbilical point if and only if $(h_{ij}(\mathcal{S}_{t_0};\bn^S)(\ou_0))=O.$
So we have the assertion (3). 
\enD
\par
We also define a family of functions
$
\widetilde{H}:U\times (LC^*\times I)\lon \R
$
by $\widetilde{H}((\ou,t),\bv)=\langle \bX(\ou,t),\widetilde{\bv}\rangle -v_0,$
where $\bv=(v_0,v_1,\dots ,v_n).$
We call $\widetilde{H}$ a {\it family of extended lightcone height functions} of $W=\bX(U\times I).$
Since $\partial\widetilde{H}/\partial u_i=\partial H/\partial u_i$ for $i=1,\dots ,s$ and
${\rm Hess} (\widetilde{h}_{(t,\sbv)})={\rm Hess} (h_{(t,\widetilde{\sbv})}),$
we have the following proposition as a corollary of Proposition 4.1.
\begin{Pro}
Let $\widetilde{H}:U\times (LC^*\times I)\lon \R$ be the extended lightcone height function of a world sheet $W=\bX(U\times I).$ Then
\par\noindent
{\rm (1)} $\widetilde{H}((\ou_0,t_0),\bv_0)=(\partial \widetilde{H}/\partial u_i)((\ou_0,t_0),\bv_0)=0\ (i=1,\dots ,s)$ if and only if
there exists a spacelike section $\bn^S$ of $N_1(\mathcal{S}_{t_0})$ such that 
$$
\bv_0=\langle \bX(\ou_0,t_0),\widetilde{\mathbb{LG}}(\mathcal{S}_{t_0};\bn^S_0)(\ou_0)\rangle \widetilde{\mathbb{LG}}(\mathcal{S}_{t_0};\bn^S_0)(\ou_0).
$$ 
\par
Suppose that $p_0=\bX(\ou_0,t_0)$, $\bv_0=\langle \bX(\ou_0,t_0),\widetilde{\mathbb{LG}}(\mathcal{S}_{t_0};\bn^S_0)(\ou_0)\rangle\widetilde{\mathbb{LG}}(
\mathcal{S}_{t_0};\bn^S_0)(\ou_0)$. Then
\par\noindent
{\rm (2)} $p_0$ is an $\bn^S_0$-parabolic point of $\mathcal{S}_{t_0}$ if and only if
${\rm det}\,{\rm Hess}(\widetilde{h}_{(t_0,\sbv_0)})(\ou_0)=0,$ 
\par\noindent
{\rm (3)} $p_0$ is a flat $\bn^S_0$-umbilical point of $\mathcal{S}_{t_0}$ if and only if
${\rm rank}\,{\rm Hess}(\widetilde{h}_{(t_0,\bm{v}_0})(\ou_0)=0$.
\end{Pro}
\demo
It follows from Proposition 4.1, (1) that
$(\partial \widetilde{H}/\partial u_i)((\ou_0,t_0),\bv_0)=0\ (i=1,\dots ,s)$ if and only if
there exists a spacelike section $\bn^S$ of $N_1(\mathcal{S}_{t_0})$ such that $\bv_0=\widetilde{\mathbb{LG}}(\mathcal{S}_{t_0};\bn^S_0)(\ou_0).$
Moreover, the condition $\widetilde{H}((\ou_0,t_0),\bv_0)=0$ is equivalent the condition 
that $v_0=\langle \bX(\ou_0,t_0),\widetilde{\mathbb{LG}}(\mathcal{S}_{t_0};\bn^S_0)(\ou_0)\rangle,$
where $\bv_0=(v_0,v_1,\dots ,v_n).$
This means that 
\[
\bv_0=\langle \bX(\ou_0,t_0),\widetilde{\mathbb{LG}}(\mathcal{S}_{t_0};\bn^S_0)(\ou_0)\rangle\widetilde{\mathbb{LG}}(\mathcal{S}_{t_0};\bn^S_0)(\ou_0).
\]
The assertions (2) and (3) directly follows from the assertions (2) and (3) of Proposition 4.1.
\enD

\par
Inspired by the above results, we define a mapping
$
\mathbb{LP}(\mathcal{S}_t):N_1(\mathcal{S}_t)\lon LC^*
$
by
\[
\mathbb{LP}(\mathcal{S}_t)((\ou,t),\bxi)=\langle \bX(\ou,t),\widetilde{\mathbb{LG}}(\mathcal{S}_{t};\bxi \rangle
\widetilde{\mathbb{LG}}(
\mathcal{S}_{t})((\ou,t),\bxi).
\]
We call it a {\it momentary lightcone pedal map} of $\mathcal{S}_t.$ Moreover, we define a map
$\mathbb{LP}:N_1(W)\lon LC^*\times I$ by
\[
\mathbb{LP}((\ou,t),\bxi)=(\mathbb{LP}(\mathcal{S}_t)((\ou,t),\bxi),t),
\]
which is called an {\it unfolded lightcone pedal map} of $W.$

\section{Contact viewpoint}
In this section we interpret the results of Propositions 4.1 and 4.2 from the view point of the contact with lightlike hyperplanes.
\par
Firstly,
we consider the relationship between the contact of a one parameter family of submanifolds with a submanifold and 
$P$-${\mathcal K}$-equivalence among functions (cf., \cite{Izudc}).  
Let $ U_i\subset \R^r$, ($i=1,2$) be open sets and $g_i:(U_i\times I, (\ou_i,t_i))\lon (\R^n,\bm{y}_i)$  immersion germs. We define $\overline{g}_i:(U_i\times I, (\ou_i,t_i))\lon (\R^n\times I,(\bm{y}_i,t_i))$
by $\overline{g}_i(\ou,t)=(g_i(\ou),t).$
We write that $(\overline{Y}_i,(\bm{y}_i,t_i))=(\overline{g}_i(U_i\times I),(\bm{y}_i,t_i)).$
Let 
$f_i:(\R^n,\bm{y}_i) \lon (\R,0)$ be submersion germs and write that $(V(f_i),\bm{y}_i)=(f_i^{-1}(0),\bm{y}_i).$
We say that {\it the contact of $\overline{Y}_1$ with the trivial family of $V(f_1)$  
at $(\bm{y}_1,t_1)$} is of the {\it same type} as {\it the contact of 
$\overline{Y}_2$ with the trivial family of $V(f_2)$ at $(\bm{y}_2,t_2)$}
if there is a diffeomorphism germ
$\Phi:(\R^n\times I,(\bm{y}_1,t_1)) \lon (\R^n\times I,(\bm{y}_2,t_2))$ of the form $\Phi (\bm{y},t)=(\phi_1(\bm{y},t),\phi _2(t))$
 such that $\Phi(\overline{Y}_1)=\overline{Y}_2$ and
$\Phi(V(f_1)\times I) = V(f_2)\times I$.
In this case we write 
$K(\overline{Y}_1,V(f_1)\times I;(\bm{y}_1,t_1)) = K(\overline{Y}_2,V(f_2)\times I;(\bm{y}_2,t_2))$.
We can show one of the parametric versions of Montaldi's theorem of contact between submanifolds as follows: 
\begin{Pro} 
We use the same notations as in the above paragraph. 
Then $K(\overline{Y}_1,V(f_1)\times I;(\bm{y}_1,t_1)) = K(\overline{Y}_2,V(f_2)\times I;(\bm{y}_2,t_2))$ 
if and only if 
$f_1 \circ g_1$ and $f_2 \circ g_2$ are $P$-${\mathcal K}$-equivalent
 {\rm (}i.e., there exists a diffeomorphism germ $\Psi :(U_1\times I,(\ou_1,t_1))\lon (U_2\times I,(\ou_2,t_2))$ of the form $\Psi (\ou,t)=(\psi_1 (\ou,t),\psi_2 (t))$ and a function germ $\lambda :(U_1\times I,(\ou_1,t_1))\lon \R$ with 
 $\lambda (\ou_1,t_1)\not= 0$ such that $(f_2\circ g_2)\circ \Phi (\ou,t) =\lambda (\ou,t)f_1\circ g_1(\ou,t)${\rm ).}
\end{Pro}
Since the proof of Proposition 5.1 is given by the arguments just along the line of the proof of the original theorem in \cite{mont1},
we omit the proof here.
\par
We now consider a function 
$
\widetilde{\mathfrak h}_{\sbv}:\R^{n+1}_1\lon {\mathbb R}
$
defined by $\widetilde{\mathfrak h}_{\sbv}(\bw )=\langle \bw ,\widetilde{\bv}\rangle -v_0,$
where $\bv=(v_0,v_1,\dots ,v_n).$
For any $\bv_0\in LC^*$, we have a lightlike hyperplane $\mathfrak{h}_{\sbv_0}^{-1}(0)=HP(\widetilde{\bv}_0 ,v_0).$
Moreover, we consider the lightlike vector $\bv_0=\mathbb{LP}(\mathcal{S}_{t_0}) ((\ou_0,t_0),\bxi_0),$
then we have 
$$
\widetilde{\mathfrak{h}}_{\sbv_0}\circ\bX (\ou_0,t_0)=\widetilde{H}(u_0,\mathbb{LP}(\mathcal{S}_{t_0}) ((\ou_0,t_0),\bxi_0)))=0.
$$
By Proposition 4.2, we also have relations that
$$
\frac{\partial \widetilde{\mathfrak{h}}_{\sbv_0}\circ\bX }{\partial u_i}(\ou_0,t_0)=\frac{\partial \widetilde{H}}{\partial u_i}((\ou_0,t_0),\mathbb{LP}(\mathcal{S}_{t_0}) ((\ou_0,t_0),\bxi_0)))=0.
$$
for $i=1,\dots ,s.$
This means that the lightlike hyperplane $\widetilde{\mathfrak{h}}_{\sbv_0}^{-1}(0)=HP(\widetilde{\bv}_0,v_0)$ is tangent to 
$\mathcal{S}_{t_0}=\bX (U\times \{t_0\})$ at $p_0=\bX (\ou_0,t_0).$
The lightlike hypersurface $HP(\widetilde{\bv}_0,v_0)$ is said to be a {\it tangent lightlike hyperplane} of $\mathcal{S}_{t_0}=\bX (U\times \{t_0\})$ at 
$p_0=\bX (\ou_0,t_0)$, which we write 
$TLP(\mathcal{S}_{t_0},\bv_0,\bxi_0)),$ where $\bv_0=\mathbb{LP}(\mathcal{S}_{t_0}) (\ou_0,t_0).$
Then we have the following simple lemma.
\begin{Lem}
Let $\bX :U\times I\lon \R^{n+1}_1$ be a world sheet. Consider two points $(p_1,\bxi_1),(p_2,\bxi_2)\in
N_1(\mathcal{S}_{t_0}) ,$ where $p_i=\bX(\ou_i,t_0)$, $(i=1,2).$
Then \[
\mathbb{LP}(\mathcal{S}_{t_0}) ((\ou_1,t_0),\bxi_1))=\mathbb{LP}(\mathcal{S}_{t_0}) ((\ou_2,t_0),\bxi_2))
\] if and only if 
$$TLP(\mathcal{S}_{t_0},\mathbb{LP}(\mathcal{S}_{t_0}) ((\ou_1,t_0),\bxi_1))=TLP(\mathcal{S}_{t_0},\mathbb{LP}(\mathcal{S}_{t_0}) ((\ou_2,t_0),\bxi_2)).$$
\end{Lem}
By definition, $\mathbb{LP}((\overline{u}_1,t_1),\bxi_1)=\mathbb{LP}((\overline{u}_2,t_2),\bxi_2)$ if and only if 
\[
t_1=t_2\ \mbox{and}\ \mathbb{LP}(\mathcal{S}_{t_1}) ((\ou_1,t_1),\bxi_1))=\mathbb{LP}(\mathcal{S}_{t_1}) ((\ou_2,t_1),\bxi_2)).
\]
Eventually, we have tools for the study of the contact between spacelike hypersurfaces and lightlike hyperplanes.
Since we have $\widetilde{h}_{\sbv}(\ou,t)=\widetilde{\mathfrak{h}}_{\sbv}\circ\bX(\ou,t),$ we have the following proposition as a
corollary of Proposition 5.1.
\begin{Pro}
Let $\bX_i : (U\times I,(\ou_i,t_i)) \lon (\R^{n+1}_1,p_i)$ $(i=1,2)$ be  world sheet germs and $\bv_i=LP(\mathcal{S}_{t_i},\mathbb{LP}(\mathcal{S}_{t_i}) ((\ou_i,t_i),\bxi_i))$
and $W_i=\bX_i(U\times I).$
Then the following conditions are equivalent:
\par\noindent
{\rm (1)} $K(\overline{W}_1, TLP(\mathcal{S}_{t_1},\bv_1,\bxi_1)\times I;(p_1,t_1))=
K(\overline{W}_2, TLP(\mathcal{S}_{t_2},\bv_2,\bxi_2)\times I;(p_2,t_2)),$ 
\par\noindent
{\rm (2)} $\widetilde{h}_{1,\sbv_1}$ and $\widetilde{h}_{2,\sbv_2}$ are $P$-$\mathcal{K}$-equivalent.\\
\end{Pro}
\section{Big wave fronts}

In this section we apply the theory of big wave fronts to the geometry of world sheets in Lorentz-Minkowski space.
Let
$
\mathcal{F}:(\R^k\times(\R^n\times\R),0)\to (\R,0)
$
be a function germ. We say that $\mathcal{F}$ is a 
{\it non-degenerate big Morse family of hypersurfaces}  
if
\[
 \Delta_*(\mathcal {F})|_{\R^k\times \R^n\times \{0\}}:
(\R^k\times \R^n\times \{0\},0)\lon (\R\times \R^k)\ \mbox{is\ non-singular},
\]
where 
\[
\Delta_*(\mathcal {F})(q,x,t)=\left(\mathcal{F}(q,x,t),\frac{\partial\mathcal{F}}{\partial q_1}(q,x,t),\dots , \frac{\partial\mathcal{F}}{\partial q_k}(q,x,t)\right).
\]
We simply say that $\mathcal{F}$ is a {\it big Morse family of hypersurfaces} if  $\Delta_*(\mathcal {F})$ is non-singular.
By definition, a non-degenerate big Morse family of hypersurfaces is a big Morse family of hypersurfaces.
Then $\Sigma_*(\mathcal{F})=\Delta (\mathcal{F})^{-1}(0)$ is
a smooth $n$-dimensional submanifold germ. 
\begin{Pro}
The extended height functions family $\widetilde{H}:U\times (LC^*\times I)\lon \R$ at any point
$(\ou_0,(\bv_0,t_0))\in \Sigma _*(\widetilde{H})$ is a non-degenerate big Morse family of hypersurfaces.
\end{Pro}
\demo
We write $\bX=(X_0,\ldots,X_n)$ and $\bv=(v_0,\ldots,v_n)\in LC^\ast$. 
	Without loss of generality, we assume that $v_n>0$. Then 
	$v_0=\sqrt{v_1^2+\cdots +v_n^2}$. 
\par
	For $\Delta^\ast \widetilde{H}=(\widetilde{H},\widetilde{H}_{u_1},\ldots,\widetilde{H}_{u_{s}})$, we prove that the map $\Delta^* \widetilde{H}|_{U\times (LC^*\times \{t_0\})}$ is submersive 
	at $(\ou_0,\bv_0,t_0)\in \Delta^\ast\widetilde{H}^{-1}(0)$. 
	Its Jacobian matrix $J\Delta^\ast \widetilde{H}|_{U\times (LC^*\times \{t_0\})}$ is 
\[
	\mbox{\Large $J\Delta^\ast \widetilde{H}|_{U\times (LC^*\times \{t_0\})}$}=\left(
		\begin{array}{c|c}
			\left(\mbox{\Large $\widetilde{H}_{u_j}$}\right)_{j=1,\ldots,s} 
			& \left(\mbox{\Large $\widetilde{H}_{v_j}$}\right)_{j=1,\ldots,n-1} 
			\\
			\hline 
			\left(\mbox{\Large $\widetilde{H}_{u_iu_j}$}\right)_{i,j=1,\ldots,s} 
			& \left(\mbox{\Large $\widetilde{H}_{u_iv_j}$}\right)_{i=1,\ldots,s,j=1,\ldots,n-1} 
		\end{array}
		\right).
		\]
We write that
\[
\mbox{\Large $B$}
=\left(
		\begin{array}{c}
			\left(\mbox{\Large $\widetilde{H}_{v_j}$}\right)_{j=1,\ldots,n-1} 
			\\
			\hline
			\left(\mbox{\Large $\widetilde{H}_{u_iv_j}$}\right)_{i=1,\ldots,s,j=1,\ldots,n-1} 
		\end{array}
		\right).
		\]
	It is enough to show that the rank of the matrix $B(\ou_0,\bv_0,t_0)$ is $s+1$. 
	By straightforward calculations, we have 
\begin{eqnarray*}
	\widetilde{H}_{v_j}(\ou,\bv,t)=-\frac{v_j}{v_0}+\frac{X_j}{v_0}-\sum_{k=1}^n\frac{v_kv_j}{v_0^3}X_k, \\ 
	\widetilde{H}_{u_iv_j}(\ou,\bv,t)=-\frac{(X_j)_{u_i}}{v_0}-\sum_{k=1}^n\frac{v_kv_j}{v_0^3}(X_k)_{u_i}, 
\end{eqnarray*}
	for $i=1,\ldots, s$ and $j=1,\ldots,n$. 
	By the condition that $\widetilde{H}(\ou_0,\bv_0,t_0)=\widetilde{H}_{u_i}(\ou_0,\bv_0,t_0)=0$ for $i$, 
	we have relations $\sum_{k=1}^n\frac{v_{0,k}}{v_{0,0}}X_k=X_0+v_{0,0}$ and 
	$\sum_{k=1}^n\frac{v_{0,k}}{v_{0,0}}(X_k)_{u_i}=(X_0)_{u_i}$ 
	where $\bv_0=(v_{0,0},\ldots,v_{0,n})$. 
	Therefore, the above formulae are
\begin{eqnarray*}
	\widetilde{H}_{v_j}(\ou_0,\bv_0,t_0)=\frac{1}{v_{0,0}}\left(X_j-2v_j-X_0\frac{v_{0,j}}{v_{0,0}}\right),\\ 
	\widetilde{H}_{u_iv_j}(\ou_0,\bv_0,t_0)=\frac{1}{v_{0,0}}\left((X_j)_{u_i}-(X_0)_{u_i}\frac{v_{0,j}}{v_{0,0}}\right), 
\end{eqnarray*}
	for $i=1,\ldots, s$ and $j=1,\ldots,n$. 

\par %
	Since $\langle \bv_0,\bv_0\rangle=\langle \bv_0,\bX_{u_i}\rangle=0$ for $i=1,\ldots, s$, 
	$\bv_0$ and $\bX_{u_i}(\ou_0,t_0)$ belong to $HP(\bv_0,0)$. 
	On the other hand, we have $\langle \bX(\ou_0,t_0)-2\bv_0+2v_{0,0}{\bf e}_0, \bv_0\rangle=-2v_{0,0}^2\neq 0$ 
	where ${\bf e}_0=(1,0,\ldots,0)$. 
	So, vectors $\bX(\ou_0,t_0)-2\bv_0+2v_{0,0}{\bf e}_0$, $\bv_0$ and $\bX_{u_i}(\ou_0,t_0)$ (for $i=1,\ldots,s$) are linearly independent. 
	Therefore the rank of following matrix 
\begin{eqnarray*}
	\mbox{\Large C}=
	\left(
		\begin{array}{c}
			\bv_0\\ X-2\bv_0+2v_{0,0}{\bf e}_0\\ X_{u_1}\\ \vdots \\ X_{u_{s}}
		\end{array}
		\right)= \left(
		\begin{array}{cccc}
			v_{0,0} & v_{0,1} & \cdots & v_{0,n} \\ 
			X_0 & X_1-2v_1 & \cdots & X_n-2v_n \\
			(X_0)_{u_1} & (X_1)_{u_1} & \cdots & (X_n)_{u_1} \\ 
			\vdots & \vdots & \ddots & \vdots \\ 
			(X_0)_{u_s} & (X_1)_{u_s} & \cdots & (X_n)_{u_s} \\ 
		\end{array}
		\right) 
\end{eqnarray*}
	is $s+2$ at $(\ou_0,\bv_0,t_0)$. 
	We subtract the first row by multiplied by $X_0/v_{0,0}$ from the second row, 
	and we also subtract the first row multiplied by $(X_0)_{u_i}/v_{0,0}$ from the $(2+i)$-th row for $i=1,\ldots,s$. 
	Then we have 
\begin{eqnarray*}
	\mbox{\Large C}'=\left(
		\begin{array}{c|c}
			v_{0,0} & v_{0,1}  \cdots v_{0,n} \\
			\hline
				\begin{array}{c}0\\ \vdots \\ 0\end{array} 
				& \mbox{\Large $B(\ou_0,\bv_0,t_0)$} \\ 
		\end{array}
		\right) 
\end{eqnarray*}
and ${\rm rank}\, C'=s+2.$
	Therefore $\mbox{rank}\, B(\ou_0,\bv_0,t_0)=s+1$.		
	 This completes the proof. 
\enD
\par
We now consider the $(n+1)$-space $\R^{n+1}=\R^n\times \R$ and coordinates of this space
are written as $(x,t)=(x_1,\dots, x_n,t)\in \R^n\times \R,$ which we distinguish space and time coordinates.
We consider the
projective cotangent 
bundle $
\pi :PT^*(\R^n\times\R)\to\R^n\times\R.
$ 
Because of the trivialization
$
PT^*(\R^n\times\R)\cong
(\R^n\times\R)\times P(({\mathbb R}^{n}\times \R)^*),
$
we have homogeneous coordinates
$$
((x_1,\dotsc ,x_n,t),[\xi _1:\cdots:\xi _n:\tau ]).
$$
Then we have the canonical contact structure $K$ on $PT^*(\R^n\times\R)$. For the definition and the basic properties of
the contact manifold $(PT^*(\R^n\times\R),K),$ see \cite[Appendix]{Izu-Pei-Sano}.
A submanifold  $i:L\subset PT^*(\R^n\times\R)$ is said to be 
a {\it big Legendrian submanifold} if $\text{dim}\, L=n$ and $di_p(T_pL)\subset K_{i(p)}$
for any $p\in L.$
We also call the map $\pi\circ i=\pi |_L: L\lon \R^n\times\R$ a {\it big Legendrian map} and the set $W(L)=\pi (L)$ a {\it big wave front} of $i:L\subset PT^*(\R^m).$
We say that a point $p\in L$ is a {\it Legendrian singular point} if ${\rm rank}\, d(\pi |_L)_p < n.$
In this case $\pi (p)$ is the singular point of $W(L).$
We call 
$$
W_t(L)=\pi_1(\pi _2^{-1}(t) \cap W(L))\quad (t\in \R)
$$
a {\it momentary front} (or, a {\it small front}) for each $t\in (\R,0),$
where $\pi_1:\R^n \times \R \to \R^n$ and $\pi_2:\R^n \times \R \to \R$ 
are the canonical projections defined by $\pi_1(x,t)=x$ and $\pi_2(x,t)=t$ respectively. 
In this sense, we call $L$ a {\it big Legendrian submanifold.}
We say that a point $p\in L$ is a {\it space-singular point} if ${\rm rank}\, d(\pi_1\circ \pi |_L)_p < n$
and a {\it time-singular point} if ${\rm rank}\, d(\pi _2\circ\pi |_L)_p=0,$ respectively.
By definition, if $p\in L$ is a Legendrian singular point, then it is a space-singular point of $L.$
\par
The {\it discriminant of the family $W_t(L)$} is defined as the image of singular points of $\pi _1|_{W(L)}.$
In the general case, the discriminant consists of three components: {\it the caustic} $C_L=\pi_1(\Sigma (W(L))$, where $\Sigma (W(L))$ is the set of singular points of $W(L)$ (i.e, the critical value set of the Legendrian mapping $\pi|_{L}$),
{\it the Maxwell stratified set} $M_L,$ the projection of self intersection points of $W(L);$
and also of the critical value set $\Delta$ of $\pi |_{W(L)\setminus \Sigma (W(L))}$ (for more detail, see \cite{Izumiya-Takahashi,Izu14,Zakalyukin95}).
We remark that $\Delta$ is not necessary the envelope of the family of smooth momentary fronts $W_t(L).$
There is a case that $\pi _2^{-1}(t)\cap W(L)$ is non-singular but $\pi _1|_{\pi _2^{-1}(t)\cap W(L)}$ has singularities,
so that $\Delta$ is the set of critical values of the family of mapping $\pi _1|_{\pi _2^{-1}(t)\cap W(L)}$ for
smooth $\pi _2^{-1}(t)\cap W(L)$. Actually, $\Delta$ is the critical value set of $\pi |_{W(L)\setminus \Sigma (W(L))}.$

\par
For any Legendrian submanifold germ $i:(L,p_0)\subset  (PT^*(\R^n\times \R),p_0),$ it is known
there exists a generating family (cf., \cite{Arnold1}). 
Let
$
\mathcal{F}:(\R^k\times(\R^n\times\R),0)\to (\R,0)
$
be a big Morse family of hypersurfaces.  Then $\Sigma_*(\mathcal{F})=\Delta (\mathcal{F})^{-1}(0)$ is
a smooth $n$-dimensional submanifold germ. 
We have a big Legendrian submanifold 
$\mathscr{L}_{\mathcal{F}}(\Sigma _*(\mathcal{F}))$ (cf., \cite{Arnold1,Zak,Zak1}), where
$$
\mathscr{L} _\mathcal{F}(q,x,t)=\left(x,t,\left[\frac{\partial \mathcal{F}}{\partial x}(q,x,t):\frac{\partial \mathcal{F}}{\partial t}(q,x,t)\right]\right),
$$
and
$$
\left[\frac{\partial \mathcal{F}}{\partial x}(q,x,t):\frac{\partial \mathcal{F}}{\partial t}(q,x,t)\right]
=\left[\frac{\partial \mathcal{F}}{\partial x_1}(q,x,t):\cdots :\frac{\partial \mathcal{F}}{\partial x_n}(q,x,t):
\frac{\partial \mathcal{F}}{\partial t}(q,x,t)\right].
$$
It is known that any big Legendrian submanifold germ can be constructed by the above method.
With this notation, the big Morse family of hypersurfaces is non-degenerate if and only if 
$(\pi_2\circ\pi)^{-1}(t)\cap\mathscr{L}_{\mathcal{F}}(\Sigma _*(\mathcal{F}))$ is a $n-1$-dimensional submanifold germ of $PT^*(\R^n\times \R)$ for
any $t\in (\R,0).$
Since $\mathscr{L}_{\mathcal{F}}(\Sigma _*(\mathcal{F}))$ is Legendrian, $(\pi_2\circ\pi)^{-1}(t)\cap\mathscr{L}_{\mathcal{F}}(\Sigma _*(\mathcal{F}))$ is an integral submanifold of the canonical contact structure $K.$
\par
We now consider an equivalence relation among big Legendrian submanifolds which preserves the discriminant of families of small fronts.
We now consider the following equivalence relation among big Legendrian submanifold germs: Let $i:(L,p_0)\subset (PT^*(\R^n\times\R),p_0)$ and $i':(L',p_0')\subset (PT^*(\R^n\times\R),p_0')$
be big Legendrian submanifold germs. 
We say that $i$ and $i'$ are {\it space-parametrized Legendrian equivalent} (or, briefly {\it $s$-$P$-Legendrian equivalent}) 
if there exist diffeomorphism germs $\Phi :(\R^n\times\R ,\overline{\pi}(p_0))\to (\R^n\times\R ,\overline{\pi}(p_0'))$
of the form $\Phi (x,t)=(\phi _1(x),\phi_2(x,t))$ such that
$\widehat{\Phi}(L)=L'$ as set germs, where 
$\widehat{\Phi}:(PT^*(\R^n\times\R),p_0) \to (PT^*(\R^n\times\R),p'_0)$ is the unique contact lift of $\Phi$.
We can also define the notion of stability of Legendrian submanifold germs with respect to $s$-$P$-Legendrian equivalence which is analogous to the stability of Lagrangian submanifold germs with respect to Lagrangian equivalence (cf. [1, Part III]).
We investigate $s$-$P$-Legendrian equivalence by using the notion of generating families of Legendrian submanifold germs.
Let $\overline{f},\overline{g}:(\R^k\times\R,0) \to (\R,0)$ be function germs.  We say that
$\overline{f}$ and $\overline{g}$ are $P$-${\cal K}$-equivalent if
there exists a diffeomorphism germ
$
\Phi :(\R^k \times \R,0)\to (\R^k \times \R,0)
$
of the form $\Phi (q,t)=(\phi _1(q,t),\phi_2(t))$ such that
$
\langle \overline{f}\circ\Phi\rangle _{{\cal E}_{k+1}}=
\langle \overline{g}\rangle _{{\cal E}_{k+1}}.
$
Let
$\mathcal{F},\mathcal{G}:(\R^k\times (\R^n\times\R),0)\to (\R,0)$ be function
germs.  We say that $\mathcal{F}$ and $\mathcal{G}$ are {\it space-$P$-${\cal K}$-equivalent} (or, briefly, {\it $s$-$P$-${\cal K}$-equivalent}) if there exists a diffeomorphism germ
$
\Psi :(\R^k\times (\R^n\times\R),0) \to (\R^k\times (\R^n\times\R),0)
$
of the form 
$
\Psi (q,x,t)=(\phi (q,x,t),\phi_1(x),\phi_2(x,t)))
$
such that
$
\langle F\circ\Psi \rangle _{{\cal E}_{k+n+1}}=
\langle G\rangle _{{\cal E}_{k+n+1}}.
$
The notion of $P$-${\cal K}$-versal deformation
plays an important role for our purpose which has been introduced in (cf.,\cite{Go-Sch,Izudc}). 
 We define the extended tangent space of $\overline{f}:(\R^k\times \R,0)\to (\R,0)$
 relative to $P$-${\cal K}$ by
$$
T_e(P\mbox{\rm -}{\cal K})(\overline{f})=\left\langle \frac{\partial \overline{f}}{\partial q_1},\dots, \frac{\partial \overline{f}}{\partial q_k},\overline{f} \right\rangle _{{\cal E}_{k+1}}+
\left\langle \frac{\partial \overline{f}}{\partial t} \right\rangle _{{\cal E}_1}.
$$
Then we say that $F$ is {\it infinitesimally $P$-${\cal K}$-versal} deformation of $\overline{f}=F|_{\R^k \times \{0\} \times \R }$ if it satisfies
$$
{\cal E}_{k+1}=T_e(P\mbox{\rm -}{\cal K})(\overline{f})+
\left\langle \frac{\partial \mathcal{F}}{\partial x_1}|_{\R^k\times\{0\}\times \R} ,\dots
,\frac{\partial \mathcal{F}}{\partial x_n}|_{\R^k\times\{0\}\times\R} \right\rangle _{\R}.
$$
We can show the following theorem analogous to those in \cite{Izu95, Zakalyukin95}.
We only remark here that the proof is analogous to the proof of \cite[Theorem in \S 21.4]{Arnold1}.
\begin{Th}
Let $\mathcal{F}:(\R^k\times(\R^n\times\R),0)\to (\R,0)$
and $\mathcal{G}:(\R^{k'}\times (\R^n\times\R),0)\to (\R,0)$
be big Morse families of hypersurfaces.  Then
\par\noindent
	{\rm ($1$)} $\mathscr{L} _\mathcal{F}(\Sigma _*(\mathcal{F}))$ and $\mathscr{L}_\mathcal{G}(\Sigma _*(\mathcal{G}))$ are $s$-$P$-Legendrian equivalent if and only if
$\mathcal{F}$ and $\mathcal{G}$ are stably $s$-$P$-${\cal K}$-equivalent.
\par\noindent
{\rm ($2$)} $\mathscr{L} _\mathcal{F}(\Sigma _*(\mathcal{F}))$ is $s$-$P$-Legendre stable if and only if $\mathcal{F}$ is an infinitesimally $P$-${\cal K}$-versal deformation of $\overline{f}=\mathcal{F}|_{\R^k\times\{0\}\times \R}.$
\end{Th}
\par
Since the Legendrian submanifold germ $i:(L,p) \subset (PT^*({\mathbb R}^n\times \R),p)$ is uniquely determined on the regular part of the big wave front $W(L),$
we have the following simple but significant property of Legendrian immersion germs \cite{Zak1}.
\begin{Pro}[Zakalyukin] Let $i:(L,p) \subset (PT^*({\mathbb R}^n\times\R),p)$, 
$i':(L',p') \subset (PT^*({\mathbb R}^n\times\R), p')$ be Legendrian submanifold germs such that regular sets of $\pi\circ i, \pi\circ i'$
 are dense respectively.
 Then $(L,p)=(L',p')$ if and only if $(W(L),\pi(p))=(W(L'),\pi (p'))$.
\end{Pro}

The assumption in Proposition 6.3 is a generic condition for $i,i'.$
Especially, if $i$ and $i'$ are $s$-$P$-Legendre stable, then these satisfy the assumption. 
Concerning the discriminant of the families of momentary fronts, we define the following equivalence relation among big wave front germs.
Let $i: (L,p_0) \subset (PT^*({\mathbb R}^n\times\R),p_0)$  and  
 $i': (L',p'_0) \subset (PT^*({\mathbb R}^n\times\R), p'_0)$  be big Legendrian submanifold germs. 
We say that $W(L)$ and $W(L')$ are {\it space-parametrized diffeomorphic} (briefly, {\it $s$-$P$-diffeomorphic}) if there exists 
a diffeomorphism germ $\Phi :(\R^n\times\R ,\overline{\pi} (p_0))\to (\R^n\times\R ,\overline{\pi} (p'_0))$
defined by $\Phi (x,t)=(\phi _1(x),\phi_2(x,t)))$ such that $\Phi (W(L))=W(L').$ 
Remark that an $s$-$P$-diffeomorphism among big wave front germs preserves the diffeomorphism types the discriminants.
By Proposition 6.3, we have the following proposition.

\begin{Pro}
Let  $i: (L,p_0) \subset (PT^*({\mathbb R}^n\times\R),p_0)$  and  
 $i': (L',p_0') \subset (PT^*({\mathbb R}^n\times\R), p_0')$  be big Legendrian submanifold germs such that regular sets of $\pi\circ i,\pi\circ i'$
 are dense respectively.
 Then $i$ and $i'$ are $s$-$P$-Legendrian equivalent  if and only if $(W(L),\pi (p_0))$ and $( W(L'),\pi (p'_0))$ are $s$-$P$-diffeomorphic.
 \end{Pro}

\begin{Rem} {\rm If we consider a diffeomorphism germ $
\Phi :(\R^n\times\R,0) \to (\R^n\times\R,0)
$
defined by 
$
\Phi (x,t)=(\phi_1(x,t), \phi _2 (t)),
$
we can define time-Legendrian equivalence among big Legendrian submanifold germs.
We can also define time-$P$-$\mathcal{K}$-equivalence among big Morse families of hypersurfaces.
By the arguments similar to the above paragraphs, we can show that these equivalence relations describe
the bifurcations of momentary fronts of big Legendrian submanifolds.
In \cite{Zak1} Zakalyukin classified generic big Legendrian submanifold germs by {\it time-Legendrian equivalence}.
The notion of time-Legendrian equivalence is a complementary notion of space-Legendrian equivalence.
}
\end{Rem}
We have the following theorem on the relation among
big Legendrian submanifolds and big wave fronts.

\begin{Th} Let $\mathcal{F}:(\R^k\times \R^n\times \R, 0)\lon (\R,0)$ and
$\mathcal{G}:(\R^{k'}\times \R^n\times \R,0)\lon (\R,0)$ be big Morse families of hypersurface  such that
$\mathscr{L}_{\mathcal{F}}(\Sigma _*(\mathcal{F}))$
and $\mathscr{L}_{\mathcal{G}}(\Sigma _*(\mathcal{G}))$ are $s$-$P$-Legendrian stable. Then the following conditions are
equivalent{\rm :}
\par\noindent
{\rm (1)} $\mathscr{L}_{\mathcal{F}}(\Sigma _*(\mathcal{F}))$
and $\mathscr{L}_{\mathcal{G}}(\Sigma _*(\mathcal{G}))$ are $s$-$P$-Legendrian equivalent,
\par\noindent
{\rm (2)} $\mathcal{F}$ and $\mathcal{G}$ are stably $s$-$P$-$\mathcal{K}$-equivalent,
\par\noindent
{\rm (3)} $\overline{f}(q,t)=\mathcal{F}(q,0,t)$ and $\overline{g}(q',t)=\mathcal{G}(q',0,t)$ are stably $P$-$\mathcal{K}$-equivalent,
\par\noindent
{\rm (4)} $W(\mathscr{L}_{\mathcal{F}}(\Sigma _*(\mathcal{F})))$ and $W(\mathscr{L}_{\mathcal{G}}(\Sigma _*(\mathcal{G})))$
are $s$-$P$-diffeomorphic.
\end{Th}
\demo
By the assertion (1) o f Theorem 6.2,  the conditions (1) and (2) are equivalent.
By definition, the condition (2) implies the condition (3).
It also follows from the definition that the condition (1) implies (4).
We remark that all these assertions hold without the assumptions of the $S$-$P$-Legendrian stability.
Generically, the condition (4) implies the condition (1) by Proposition 6.4.
Of course, it holds under the assumption of $S$-$P$-Legendrian stability.
By the assumption of $s$-$P$-Legendrian stability, the big Morse families of
hypersurface $\mathcal{F}$ and $\mathcal{G}$ are infinitesimally $P$-$\mathcal{K}$-versal deformations
of  $\overline{f}$ and  $\overline{g}$, respectively (cf., Theorem 6.2, (2)).
By the uniqueness result of the infinitesimally $P$-$\mathcal{K}$-versal deformations (cf., \cite{Izudc}), the condition (3) implies the condition (2). This completes the proof.
\enD
\begin{Rem}{\rm 
(1) If $k=k'$ and $q=q'$ in the above theorem, we can remove the word \lq\lq stably\rq\rq
\ in the conditions (2),(3).
\par\noindent
(2) $s$-$P$-Legendrian stability for $\mathscr{L}_{\mathcal{F}}(\Sigma _*(\mathcal{F}))$ is generic for $n\leq 5.$ 
\par\noindent
(3) By the remark in the proof of the above theorem, the conditions (1) and (4) are equivalent
generically for a general dimension $n$ without the assumption on $s$-$P$-Legendrian stability.
Therefore, the conditions (1),(2) and (4) are all equivalent to each other generically.
}
\end{Rem}

\par
We now return to our situation. Since the extended lightcone height functions family $\widetilde{H}:U\times (LC^*\times I)\lon \R$
is a non-degenerate big Morse family of hypersurfaces, we have the corresponding big Legendrian submanifold
$\mathscr{L}_{\widetilde{H}}(\Sigma _*(\widetilde{H}))\subset PT^*(LC^*\times I).$
By Proposition 4.2, we have
\begin{eqnarray*}
\Sigma _*(\widetilde{H}))&=&\{(\ou,\mathbb{LP}(\mathcal{S}_t)((\ou,t),\bxi),t)\ |\  (\ou,t)\in U\times I, \bxi\in N_1[\mathcal{S}_t]_p,
p=\bX(\ou,t) \} \\
&=& \{(u,\mathbb{LP}((\ou,t),\bxi))\ |\ (\ou,t)\in U\times I, \bxi\in N_1[\mathcal{S}_t]_p,
p=\bX(u,t) \}.
\end{eqnarray*}
It follows that $\pi (\mathcal{L}_{\widetilde{H}}(\Sigma _*(\widetilde{H})))=\mathbb{LP}(N_1(W))\subset LC^*\times I.$
Therefore, the image of the unfolded lightcone pedal is a big wave front. 

\par
We apply the above theorem to our situation.
\begin{Th}
Let $\bX_i : (U\times I,(\ou_i,t_i)) \lon (\R^{n+1}_1,p_i)$ $(i=1,2)$ be  world sheet germs and $\bv_i=\mathbb{LP}(\mathcal{S}_{t_i}) ((\ou_i,t_i),\bxi_i))$
and $W_i=\bX_i(U\times I).$
Suppose that the Legendrian submanifold germs $\mathcal{L}_{\widetilde{H}_i}(\Sigma _*(\widetilde{H}_i))\subset PT^*(LC^*\times I)$
 are $s$-$P$-Legendrian stable.
Then the following conditions are equivalent:
\par\noindent
{\rm (1)} $\mathcal{L}_{\widetilde{H}_1}(\Sigma _*(\widetilde{H}_1))$ and $\mathcal{L}_{\widetilde{H}_2}(\Sigma _*(\widetilde{H}_2))$ are $s$-$P$-Legendrian equivalent,
\par\noindent
{\rm (2)} $\widetilde{h}_{1,\sbv_1}$ and $\widetilde{h}_{2,\sbv_2}$ are $P$-$\mathcal{K}$-equivalent,
\par\noindent
{\rm (3)} $\widetilde{H}_1$ and $\widetilde{H}_2$ are $s$-$P$-$\mathcal{K}$-equivalent,
\par\noindent
{\rm(4)} $\mathbb{LP}_1(N_1(W_1))$ and $\mathbb{LP}_2(N_1(W_2))$ are $s$-$P$-diffeomorphic,
\par\noindent
{\rm (5)} $K(\overline{W}_1, TLP(\mathcal{S}_{t_1},\bv_1,\bxi_1)\times I;(p_1,t_1))=
K(\overline{W}_2, TLP(\mathcal{S}_{t_1},\bv_2,\bxi_2)\times I;(p_2,t_2)).$
\end{Th}
\demo
Since $\mathbb{LP}_i(N_1(W_i))$ are big wave fronts of $\mathcal{L}_{\widetilde{H}_i}(\Sigma _*(\widetilde{H}_i))$ $(i=1,2)$
respectively, we can apply Theorem 6.6 and obtain that the conditions (1), (2), (3) and (4) are equivalent.
By Proposition 5.3, the conditions (2) and (5) are equivalent.
This completes the proof.
\enD

\end{document}